\documentclass[11pt]{amsart}

\usepackage{amsmath,amsthm}
\usepackage{amssymb}
\usepackage{euscript}
\usepackage[frame,ps,matrix,arrow,curve,rotate]{xy}

\numberwithin{equation}{subsection}  

\newcommand{\sqsp}{\renewcommand{\baselinestretch}{1.1}\tiny\normalsize}

\newtheorem{thm}[subsection]{Theorem}
\newtheorem{lemma}[subsection]{Lemma}
\newtheorem{prop}[subsection]{Proposition}
\newtheorem{cor}[subsection]{Corollary}

\newtheorem*{thm-a}{Main Theorem}
\newtheorem*{cor-a}{Corollary}

\theoremstyle{definition}
\newtheorem{definition}[subsection]{Definition}

\newcommand{\cat}[1]{{\EuScript #1}}

\newcommand{\cE}{\cat{E}}
\newcommand{\cF}{\cat{F}}

\newcommand{\bR}{\mathbf{R}}

\newcommand{\bZ}{\mathbf{Z}}

\DeclareMathOperator{\Id}{Id}
\DeclareMathOperator{\End}{End}
\DeclareMathOperator{\Hom}{Hom}
\DeclareMathOperator{\Der}{Der}
\DeclareMathOperator{\Diff}{Diff}
\DeclareMathOperator{\Set}{\textbf{Set}}
\DeclareMathOperator{\Ob}{Ob}

\newcommand{\biglbrack}{\biggl \lbrack}  
\newcommand{\bigrbrack}{\biggr \rbrack}  

\begin{document}
\title{Deformation of algebras over the Landweber-Novikov algebra}
\author{Donald Yau}

\begin{abstract}
An algebraic deformation theory of algebras over the Landweber-Novikov algebra is obtained.  \end{abstract}

\email{dyau@math.ohio-state.edu}
\address{Department of Mathematics, The Ohio State University Newark, 1179 University Drive, Newark, OH 43055}

\maketitle
\sqsp


\section{Introduction}
\label{sec:intro}

The Landweber-Novikov algebra $S$ was originally studied by Landweber \cite{landweber} and Novikov \cite{novikov} as a certain algebra of stable cohomology operations on complex cobordism $MU^*(-)$.   In fact, it is known that every stable cobordism operation can be written uniquely as an $MU^*$-linear combination in the Landweber-Novikov operations.  The Landweber-Novikov operations act stably, and hence additively, on the cobordism $MU^*(X)$ of a space $X$.  However, they are, in general, not multiplicative on $MU^*(X)$.  Instead, their actions on a product of two elements in $MU^*(X)$ satisfy the so-called Cartan formula, analogous to the formula of the same name in ordinary mod $2$ cohomology.  This structure - the $S$-module structure together with the Cartan formula on products - makes $MU^*(X)$ into what is called an \emph{algebra over the Landweber-Novikov algebra}, or an $S$-algebra for short.  These $S$-algebras are therefore of great importance in algebraic topology.

The algebra $S$ has also appeared in other settings.  For example,  Bukhshtaber and Shokurov \cite{bs} showed that the Landweber-Novikov algebra is isomorphic to the algebra of left invariant differentials on the group $\Diff_1(\bZ)$.  Denoting the group of formal diffeomorphisms on the real line by $\Diff_1(\bR)$, the group $\Diff_1(\bZ)$ is the subgroup generated by the formal diffeomorphisms with integer coefficients.  This theme that the Landweber-Novikov algebra is an ``operation algebra" is echoed in Wood's paper \cite{wood}.  Wood constructed $S$ as a certain algebra of differential operators on the integral polynomial ring on countably infinitely many variables.  There are even connections between the Landweber-Novikov algebra and physics, as the work of Morava \cite{morava} demonstrates.

The purpose of this paper is to study algebras over the Landweber-Novikov algebra from the specific view point of algebraic deformations.  We would like to deform an $S$-algebra $A$ with respect to the Landweber-Novikov operations on $A$, keeping the algebra structure on $A$ unaltered.  The resulting deformation theory is described in cohomological and obstruction theoretic terms.

The original theory of deformations of associative algebras was developed by Gerstenhaber in a series of papers \cite{ger1,ger2,ger3,ger4}.  It has since been extended in many different directions, and many kinds of algebras now have their own deformation theories.

The Landweber-Novikov algebra is actually a Hopf algebra, as the referee pointed out.  Therefore, what we are considering in this paper is really an instance of deformation of a module algebra over a Hopf algebra.  It would be nice to extend the results in the current paper to this more general setting.

A description of the rest of the paper follows.

The following section is preliminary in nature.  We recall the Landweber-Novikov algebra $S$ and algebras over it.  Our deformation theory depends on the cohomology of a certain cochain complex $\cF^*$, which is constructed in the next section as well.

In Section \ref{sec:formal} we introduce the notions of a formal deformation and of a formal automorphism.  The latter is used to defined equivalence of formal deformations.  The main point of that section is Theorem \ref{thm:inf}, which identifies the ``infinitesimal" of a formal deformation with an appropriate cohomology class in $H^1(\cF^*)$.  Intuitively, the ``infinitesimal" is the initial velocity of the formal deformation.

Section \ref{sec:extending} begins with a discussion of formal automorphisms of finite order and how such objects can be extended to higher order ones.  See Theorem \ref{thm:ext} and Corollary \ref{cor:ext}.  These results are needed to study rigidity.  An $S$-algebra $A$ is called \emph{rigid} if every formal deformation of $A$ is equivalent to the trivial one.  The main result there is Corollary \ref{cor2:rigid}, which states that an $S$-algebra $A$ is rigid, provided that both $H^1(\cF^*)$ and $HH^2(A)$ are trivial.  Here $HH^2(A)$ denotes the second Hochschild cohomology of $A$, as an algebra over the ring of integers, with coefficients in $A$ itself.

In Section \ref{sec:integration}, we identify the obstructions to extending a $1$-cocycle in $\cF^*$ to a formal deformation.  This is done by considering formal deformations of finite orders and identifying the obstructions to extending such objects to higher order ones.  This obstruction turns out to be in $H^2(\cF^*)$; see Theorem \ref{thm1:int}.  As a result, the vanishing of $H^2(\cF^*)$ implies that every $1$-cocycle occurs as the infinitesimal of a formal deformation (Corollary \ref{cor2:int}).  The paper ends with Theorem \ref{thm2:int}, which shows that the vanishing of a certain class in $H^1(\cF^*)$ implies that two order $m + 1$ extensions of an order $m$ formal deformation are equivalent.


\section{The Landweber-Novikov algebra and the complex $\cF^*$}
\label{sec:prelim}

The purpose of this preliminary section is to recall the Landweber-Novikov algebra $S$ and the notion of an algebra over it.  Then we construct a cochain complex $\cF^*$ associated to an algebra over the Landweber-Novikov algebra.  This complex will be used in later sections to study algebraic deformations of $S$-algebras.

\subsection{The Landweber-Novikov algebra}
\label{subsec:ln algebra}

References for this subsection are \cite{landweber,novikov}, where the Landweber-Novikov algebra was first introduced.  Wood's paper \cite{wood} has an alternative description of it as a certain algebra of differential operators.  The book \cite[Ch.\ I]{adams} by Adams is also a good reference.

The Landweber-Novikov algebra $S$ is generated by certain elements $s_\alpha$, the Landweber-Novikov operations, indexed by the exponential sequences, which we first recall.

An \emph{exponential sequence} is a sequence
   \[
   \alpha ~=~ (\alpha_1, \alpha_2, \ldots \,)
   \]
of non-negative integers in which all but a finite number of the $\alpha_i$ are $0$.  When $\alpha$ and $\beta$ are two exponential sequences, their sum $\alpha + \beta$ is defined componentwise.  Denote by $\cE$ the set of all exponential sequences.

For each exponential sequence $\alpha$, there is a stable cobordism cohomology operation
   \[
   s_\alpha \colon MU^*(-) \to MU^{*^\prime}(-),
   \]
called a Landweber-Novikov operation.  The composition $s_\alpha s_\beta$ of any two Landweber-Novikov operations satisfies the \emph{product formula}, which expresses it uniquely as a finite $\bZ$-linear combination,
   \begin{equation}
   \label{eq:product}
   s_\alpha s_\beta ~=~ \sum_{\gamma \in P(\alpha, \beta)} \, n_\gamma s_\gamma.
   \end{equation}
See, for example, Adams \cite[\S\S 5,6]{adams} or Wood \cite[Thm.\ 4.2]{wood} for a proof of this fact.  For our purposes, we do not really need to know what the integers $n_\gamma$ are, as long as they satisfy the obvious associativity condition.

Thus, using the product formula \eqref{eq:product}, the free abelian group
   \begin{equation}
   \label{eq:landweber novikov algebra}
   S ~=~ \bigoplus_{\alpha \in \cE} \, \bZ s_\alpha
   \end{equation}
can be equipped with an algebra structure with composition as the product.    This is what we referred to as the \emph{Landweber-Novikov algebra}.

The cobordism $MU^*(X)$ of a space or a spectrum $X$ is automatically an $S$-module.  In fact, it is known that every stable cobordism cohomology operation can be written uniquely as an $MU^*$-linear combination of the Landweber-Novikov operations.  (Here $MU^*$ denotes $MU^*(\mathrm{pt})$ as usual.)  If $X$ is a space, then $MU^*(X)$ is not just a group but an algebra as well.  On a product of two elements, the Landweber-Novikov operations satisfy the \emph{Cartan formula},
   \begin{equation}
   \label{eq:cartan}
   s_\alpha(ab) ~=~ \sum_{\beta + \gamma = \alpha}\, s_\beta(a)s_\gamma(b).
   \end{equation}
Here $\alpha \in \cE$ and $a, b \in MU^*(X)$.  We, therefore, make the following definition.

\begin{definition}
\label{def:S-alg}
By an \emph{algebra over the Landweber-Novikov algebra}, or an \emph{$S$-algebra} for short, we mean a commutative ring $A$ with $1$ that comes equipped with an $S$-module structure such that the Cartan formula \eqref{eq:cartan} is satisfied for all $\alpha \in \cE$ and $a, b \in A$.
\end{definition}

Given a commutative ring $A$, let $\End(A)$ denote the algebra of additive self-maps of $A$, where product is composition of self maps.  If $f$ and $g$ are in $\End(A)$, then $fg$ always means the composition $f \circ g$.   An $S$-algebra structure on $A$ is equivalent to a function
  \begin{equation}
  \label{eq:s}
  s_* \colon \cE ~\to~ \End(A),
  \end{equation}
assigning to each exponential sequence $\alpha$ an additive operation $s_\alpha$ on $A$, which satisfies the product formula \eqref{eq:product} on compositions and the Cartan formula \eqref{eq:cartan} on products in $A$.

\subsection{The complex $\cF^*$}
\label{subsec:F}

The algebraic deformation theory of $S$-algebras discussed in later sections depends on a certain cochain complex $\cF^*$, which we now construct.

First we need some notations.  Let $A$ be a commutative ring.  Denote by $\Der(A)$ the abelian group of derivations on $A$.  Recall that a \emph{derivation} on $A$ is an additive self-map $f \colon A \to A$ which satisfies the condition
   \[
   f(ab) ~=~ af(b) ~+~ f(a)b
   \]
for all $a, b \in A$.  For a positive integer $n$, we use the shorter notation $A^{\otimes n}$ to denote the tensor product over the ring of integers of $A$ with itself $n$ times.  The group of additive maps $A^{\otimes n} \to A$ is written $\Hom(A^{\otimes n}, A)$.  When $n = 1$, we also write $\End(A)$ for $\Hom(A, A)$.  Each set $\Hom(A^{\otimes n}, A)$ has a natural abelian group structure.  Namely, given $f, g \colon A^{\otimes n} \to A$ and $b \in A^{\otimes n}$, the sum $(f + g)$ sends $b$ to $f(b) + g(b)$.

Let $\Set$ be the category of sets.  Given two sets $C$ and $D$, the set of functions from $C$ to $D$ is written $\Set(C, D)$.  If $D$ is an abelian group, then, just as in the previous paragraph, the set $\Set(C, D)$ is naturally equipped with an abelian group structure as well.

Let $A$ be an algebra over the Landweber-Novikov algebra $S$ (see Definition \ref{def:S-alg}).  Consider the Landweber-Novikov operations on $A$ as a function $s_*$ as in \eqref{eq:s}.  Recall that $\cE$ is the set of all exponential sequences.  Also recall the sets $P(\alpha, \beta)$ from the product formula \eqref{eq:product}.    Denote by $\cE^n$ the Cartesian product $\cE \times \cdots \times \cE$ ($n$ factors).  Just as in the case of $s_*$, if $f$ is a function with domain $\cE$, we will write $f_\alpha$ instead of $f(\alpha)$ for $\alpha \in \cE$.  We are now ready to define the cochain complex $\cF^* = \cF^*(A)$ of abelian groups.

We make the following definitions.
\begin{itemize}
\item $\cF^0(A) = \Der(A)$.
\item $\cF^1(A) = \Set(\cE, \End(A))$.
\item For integers $n \geq 2$, set 
   \[
   \cF^n(A) \,=\, \cF^n_0(A) \,\times\, \cF^n_1(A), 
   \]
where
   \[
   \begin{split}
   \cF^n_0(A) &\,=\, \Set(\cE^n, \End(A)), \\
   \cF^n_1(A) &\,=\, \Set(\cE, \Hom(A^{\otimes n}, A)).
   \end{split}
   \]
\end{itemize}
Now we define the differentials.

\begin{itemize}
\item $d^0(\varphi) = s_* \varphi - \varphi s_*$ for $\varphi \in \cF^0(A)$.
\item For $n \geq 1$, $d^n = (d^n_0, d^n_1)$, where
   \[
   \begin{split}
   d^n_0 & \colon \cF^n_0(A) \,\to\, \cF^{n+1}_0(A), \\
   d^n_1 & \colon \cF^n_1(A) \,\to\, \cF^{n+1}_1(A)
   \end{split}
   \]
are defined as follows.  (Here we have $\cF^1_0(A) = \cF^1_1(A) = \cF^1(A)$.)  Suppose that $f = (f_0, f_1)$ is an element of $\cF^n(A)$ for some $n \geq 2$ (or just $f$ when $n = 1$), $\mathbf{x} = (\alpha_1, \ldots, \alpha_{n+1}) \in \cE^{n+1}$, $\alpha \in \cE$, and $\mathbf{a} = a_1 \otimes \cdots \otimes a_{n+1} \in A^{\otimes (n+1)}$.  Then we set
   \begin{equation}
   \label{eq:dn0}
   \begin{split}
   (d^n_0 f_0)(\mathbf{x}) \,=\, s_{\alpha_1} & f_0(\alpha_2, \ldots, \alpha_{n+1}) \\
   & ~+~ \sum_{i=1}^n (-1)^i\biglbrack\sum_{\beta \in P(\alpha_i,\alpha_{i+1})} n_\beta f_0(\cdots, \alpha_{i-1}, \beta, \alpha_{i+2}, \ldots)\bigrbrack \\
   & ~+~ (-1)^{n+1} f_0(\alpha_1, \ldots, \alpha_n)s_{\alpha_{n+1}},
   \end{split}
   \end{equation}
where $P(-, -)$ is as in the product formula \eqref{eq:product}, and
   \begin{equation}
   \label{eq:dn1}
   \begin{split}
   (d^n_1 f_1)(\alpha)(\mathbf{a}) \,=\, &\sum_{\beta + \gamma = \alpha}  s_\beta(a_1)f_1(\gamma)(a_2 \otimes \cdots \otimes a_{n+1}) \\
   & ~+~ \sum_{i=1}^n (-1)^i f_1(\alpha)(\cdots \otimes a_{i-1} \otimes (a_ia_{i+1}) \otimes a_{i+2} \cdots) \\
   & ~+~ (-1)^{n+1} \sum_{\beta + \gamma = \alpha} f_1(\beta)(a_1 \otimes \cdots \otimes a_n)s_\gamma(a_{n+1}).
   \end{split}
   \end{equation}
In these definitions, when $n = 1$, both $f_0$ and $f_1$ are interpreted as $f$.
\end{itemize}

\begin{prop}
\label{prop:F cochain}
($\cF^*, d^*)$ is a cochain complex of abelian groups.
\end{prop}

\begin{proof}
It is straightforward to check that $d^1d^0 = 0$ by direct inspection.

For $i = 0, 1$ and $n \geq 1$, the formulas above allow one to rewrite $d^n_i$ as an alternating sum
$d^n_i \,=\, \sum_{j=0}^{n+1} (-1)^j \partial^n_i \lbrack j \rbrack,$ corresponding to the $n + 2$ terms in the respective formulas.  It is straightforward to check that, for $n \geq 1$, the cosimplicial identities, 
   \[
   \partial^{n+1}_i \lbrack l \rbrack \circ  \partial^n_i \lbrack k \rbrack 
   \,=\, 
   \partial^{n+1}_i\lbrack k \rbrack \circ \partial^n_i \lbrack l - 1 \rbrack \quad (0 \leq k < l \leq n + 2)
   \]
hold.  Therefore, it follows that $\cF^*(A)$ is a cochain complex.
\end{proof}

From now on, whenever we speak of ``cochains," ``cocycles," ``coboundaries," and ``cohomology classes," we are referring to the cochain complex $\cF^* = (\cF^*, d^*)$, unless stated otherwise.  The $i$th cohomology group of $\cF^*$ will be denoted by $H^i(\cF^*)$.


\section{Formal deformations and automorphisms}
\label{sec:formal}

The purposes of this section are (1) to introduce formal deformations and automorphisms of an $S$-algebra and (2) to identify the ``infinitesimal" of a formal deformation with an appropriate cohomology class in the complex $\cF^*$.

\subsection{Formal deformations}
\label{subsec:formal def}

Throughout this section, let $A$ be an arbitrary but fixed $S$-algebra.  Recall that the $S$-algebra structure on $A$ can be characterized as a function $s_*$ as in \eqref{eq:s} from the set $\cE$ of exponential sequences to the algebra of additive self maps of $A$.  In addition, this function satisfies the Cartan formula \eqref{eq:cartan} and the product formula \eqref{eq:product}.  We will deform the Landweber-Novikov operations on $A$ with respect to these two properties.

We define a \emph{formal deformation} of $A$ to be a formal power series in the indeterminate $t$,
   \begin{equation}
   \label{eq:def}
   \sigma^t_* ~=~ s_* ~+~ ts^1_* ~+~ t^2s^2_* ~+~ t^3s^3_* ~+~ \cdots,
   \end{equation}
in which each $s^i_* \in \cF^1(A)$ $(i \geq 1)$, i.e.\ is a function $\cE \to \End(A)$, satisfying the following two properties (with $s^0_* = s_*$ and $s^i_*(\alpha) = s^i_\alpha$):
   \begin{itemize}
   \item The Cartan formula: For every $\alpha \in \cE$ and $a, b \in A$, the equality
   \begin{equation}
   \label{eq:Cartan}
   \sigma^t_\alpha(ab) ~=~ \sum_{\beta + \gamma = \alpha}\sigma^t_\beta(a)\sigma^t_\gamma(b)
   \end{equation}
   of power series holds.
   \item The Product formula: For $\alpha,\beta \in \cE$, the equality
   \begin{equation}
   \label{eq:Product}
   \sigma^t_\alpha \sigma^t_\beta ~=~ \sum_{\gamma \in P(\alpha, \beta)} n_\gamma \sigma^t_\gamma
   \end{equation}
   of power series holds.
   \end{itemize}

We pause to make a few remarks.  First, the superscript $i$ in $s^i_*$ is an index, not an exponent, whereas $t^i$ is the $i$th power of the indeterminate $t$.  Second, sums and products of two power series are taken in the usual way, with $t$ commuting with every term in sight.  The coefficients in \eqref{eq:Cartan} and \eqref{eq:Product} are in $A$ and $\End(A)$, respectively, and their calculations are done in the appropriate rings.  In particular, multiplying out the right-hand side of the equation, one observes that the Cartan formula \eqref{eq:Cartan} is equivalent to the equality
  \begin{equation}
  \label{eq:Cartan'}
  s^n_\alpha(ab) ~=~ \sum_{i=0}^n \sum_{\beta + \gamma = \alpha} s^i_\beta(a)s^{n-i}_\gamma(b)
  \end{equation}
in $A$ for all $n \geq 0$, $\alpha \in \cE$, and $a, b \in A$.  Similarly, unwrapping the left-hand side of the equation, one observes that the Product formula \eqref{eq:Product} is equivalent to the equality
   \begin{equation}
   \label{eq:Product'}
   \sum_{i=0}^n s^i_\alpha s^{n-i}_\beta ~=~ \sum_{\gamma \in P(\alpha, \beta)} n_\gamma s^n_\gamma
   \end{equation}
in $\End(A)$ for all $n \geq 0$ and $\alpha, \beta \in \cE$.  When $n = 0$, \eqref{eq:Cartan'} and \eqref{eq:Product'} are just the original Cartan formula \eqref{eq:cartan} and product formula \eqref{eq:product}, respectively, for $s_* = s^0_*$.

Setting $t = 0$ in the formal deformation $\sigma^t_*$, we obtain $\sigma^0_* = s_*$.  So we can think of $\sigma^t_*$ as a one-parameter curve with $s_*$ at the original.  We, therefore, call $s^1_*$ the \emph{infinitesimal}, since it is the ``initial velocity" of the formal deformation $\sigma^t_*$.

\subsection{Formal automorphisms}
\label{subsec:formal aut}

In order to identify the infinitesimal as an appropriate cohomology class in $\cF^*$, we need a notion of equivalence of formal deformations.

By a \emph{formal automorphism} on $A$, we mean a formal power series
   \begin{equation}
   \label{eq:aut}
   \Phi_t ~=~ 1 ~+~ t\phi_1 ~+~ t^2\phi_2 ~+~ t^3\phi_3 ~+~ \cdots,
   \end{equation}
where $1 = \Id_A$ and each $\phi_i \in \End(A)$, satisfying \emph{multiplicativity},
   \begin{equation}
   \label{eq:mult}
   \Phi_t(ab) ~=~ \Phi_t(a)\Phi_t(b)
   \end{equation}
for all $a, b \in A$.

The same rules of dealing with power series apply here as well.  In particular, multiplicativity is equivalent to the equality
   \begin{equation}
   \label{eq:mult'}
   \phi_n(ab) ~=~ \sum_{i=0}^n \phi_i(a)\phi_{n-i}(b)
   \end{equation}
for all $n \geq 0$ and $a, b \in A$, in which $\phi_0 = 1$.  The condition when $n = 0$ is trivial, as it only says that the identity map on $A$ is multiplicative.  When $n = 1$, the condition is
   \begin{equation}
   \label{eq:phi1}
   \phi_1(ab) ~=~ a\phi_1(b) ~+~ \phi_1(a)b,
   \end{equation}
which is equivalent to say that $\phi_1$ is a derivation on $A$.  More generally, if $\phi_1 = \phi_2 = \cdots = \phi_k = 0$, then $\phi_{k+1}$ is a derivation on $A$.

It is an easy exercise in induction to see that a formal automorphism $\Phi_t$ has a unique formal inverse
   \begin{equation}
   \label{eq:inv}
   \Phi^{-1}_t ~=~ 1 ~-~ t\phi_1 ~+~ t^2(\phi_1^2 - \phi_2) ~+~ t^3(-\phi_1^3 + \phi_1\phi_2 + \phi_2\phi_1 - \phi_3) ~+~ \cdots,
   \end{equation}
for which
   \[
   \Phi_t\Phi^{-1}_t ~=~ 1 ~=~ \Phi^{-1}_t\Phi_t.
   \]
The coefficient of $t^n$ in $\Phi^{-1}_t$ is an integral polynomial in $\phi_1, \cdots, \phi_n$.  Moreover, the multiplicativity of $\Phi_t$ implies that of $\Phi^{-1}_t$.  Indeed, for elements $a, b \in A$, we have
   \[
   ab ~=~ (\Phi_t \Phi^{-1}_t(a))(\Phi_t \Phi^{-1}_t(b)) ~=~ \Phi_t(\Phi^{-1}_t(a) \Phi^{-1}_t(b)),
   \]
which implies that $\Phi^{-1}_t$ is multiplicative.

We record these facts as follows.

\begin{lemma}
\label{lem1:aut}
Let $\Phi_t = 1 + t\phi_1 + t^2\phi_2 + \cdots$ be a formal automorphism on $A$.  Then the first non-zero $\phi_i$ $(i \geq 1)$ is a derivation on $A$.  Moreover, the formal inverse $\Phi^{-1}_t$ of $\Phi_t$ is also a formal automorphism on $A$.
\end{lemma}

Now if $s^n_*$ is a function $\cE \to \End(A)$, i.e.\ a $1$-cochain, and if $f$ and $g$ are in $\End(A)$, then we have a new $1$-cochain $fs^n_* g$ with
   \[
   (fs^n_*g)(\alpha) ~=~ fs^n_\alpha g
   \]
in $\End(A)$ for $\alpha \in \cE$.  Therefore, it makes sense to consider the formal power series $\Phi^{-1}_t \sigma^t_* \Phi_t$ whenever $\Phi_t$ is a formal automorphism.

\begin{prop}
\label{prop:aut}
Let $\sigma^t_*$ and $\Phi_t$ be, respectively, a formal deformation and a formal automorphism of $A$.  Then the formal power series $\Phi^{-1}_t \sigma^t_* \Phi_t$ is also a formal deformation of $A$.
\end{prop}

\begin{proof}
We need to check the Cartan formula \eqref{eq:Cartan} and the Product formula \eqref{eq:Product} for $\tilde{\sigma}^t_* = \Phi^{-1}_t \sigma^t_* \Phi_t$.  For the Cartan formula,  we have
   \[
   \begin{split}
   \tilde{\sigma}^t_\alpha(ab)
   &~=~ (\Phi^{-1}_t \sigma^t_\alpha)(\Phi_t(a) \Phi_t(b)) \\
   &~=~ \Phi^{-1}_t\biggl(\sum_{\beta + \gamma = \alpha}(\sigma^t_\beta\Phi_t(a))(\sigma^t_\gamma \Phi_t(b))\biggr) \\
   &~=~ \sum_{\beta + \gamma = \alpha} (\Phi^{-1}_t\sigma^t_\beta\Phi_t(a))(\Phi^{-1}_t\sigma^t_\gamma\Phi_t(b)).
   \end{split}
   \]
We have used the Cartan formula for $\sigma^t_*$ and the multiplicativity of both $\Phi_t$ and $\Phi^{-1}_t$, extended to power series.  The Product formula is equally easy to verify.
\end{proof}

Given two formal deformations $\sigma^t_*$ and $\tilde{\sigma}^t_*$ of $A$, say that they are \emph{equivalent} if and only if there exists a formal automorphism $\Phi_t$ on $A$ such that
   \begin{equation}
   \label{eq:equivalence}
   \tilde{\sigma}^t_* ~=~ \Phi^{-1}_t \sigma^t_* \Phi_t.
   \end{equation}
By Lemma \ref{lem1:aut} this is a well-defined equivalence relation on the set of formal deformations of $A$.

Here is the main result of this section, which identifies the infinitesimal with an appropriate cohomology class in $H^1(\cF^*)$.  Recall that cocycles, coboundaries, cochains, and cohomology classes are all taken in the cochain complex $\cF^* = \cF^*(A)$.

\begin{thm}
\label{thm:inf}
Let $\sigma^t_* = \sum_n t^n s^n_*$ be a formal deformation of $A$.  Then the infinitesimal $s^1_*$ is a $1$-cocycle, i.e.\ $d^1s^1_* = 0$.  Moreover, the cohomology class $\lbrack s^1_* \rbrack$ is an invariant of the equivalence class of $\sigma^t_*$.

More generally, if $s^1_* = \cdots = s^k_* = 0$ for some positive integer $k$, then  $s^{k+1}_*$ is a $1$-cocycle.
\end{thm}

\begin{proof}
To show that $s^1_*$ is a $1$-cocycle, we need to prove that $d^1_i s^1_* = 0$ for $i = 0, 1$.  When $n = 1$ the Product formula \eqref{eq:Product'} states that
   \[
   s_\alpha s^1_\beta ~+~ s^1_\alpha s_\beta ~=~ \sum_{\gamma \in P(\alpha, \beta)}n_\gamma s^1_\gamma,
   \]
and so
   \[
   (d^1_0s^1_*)(\alpha, \beta) ~=~ s_\alpha s^1_\beta ~-~ \sum_{\gamma \in P(\alpha, \beta)}n_\gamma s^1_\gamma ~+~ s^1_\alpha s_\beta ~=~ 0.
   \]
Similarly, the Cartan formula \eqref{eq:Cartan'} when $n = 1$ states that
   \[
   s^1_\alpha(ab) ~=~ \sum_{\beta + \gamma = \alpha} (s_\beta(a)s^1_\gamma(b) ~+~ s^1_\beta(a)s_\gamma(b)),
   \]
which implies that
   \[
   (d^1_1s^1_*)_\alpha(a \otimes b) ~=~ 0,
   \]
as desired.  This shows that $s^1_*$ is a $1$-cocycle.

If $s^1_* = \cdots = s^k_* = 0$, then the same argument as above shows that $s^{k+1}_* = 0$.

Now suppose that $\tilde{\sigma}^t_* = \sum_n t^n \tilde{s}^n_*$ is a formal deformation of $A$ that is equivalent to $\sigma^t_*$.  This means that there exists a formal automorphism $\Phi_t$ such that
   \[
   \begin{split}
   \tilde{\sigma}^t_*
   &~=~ \Phi^{-1}_t \sigma^t_* \Phi_t \\
   &~\equiv~ s_* ~+~ t(s^1_* + s_*\phi_1 - \phi_1s_*) \qquad \pmod{t^2} \\
   &~\equiv~ s_* ~+~ t(s^1_* + d^0\phi_1) \qquad \pmod{t^2}.
   \end{split}
   \]
In particular, the $1$-cocycle $(\tilde{s}^1_* - s^1_*) \in \cF^1$ is a $1$-coboundary $d^0\phi_1$.  (Remember that $\phi_1$ is a derivation on $A$, which is therefore a $0$-cochain.)  Thus, the cohomology classes, $\lbrack \tilde{s}^1_* \rbrack$ and $\lbrack s^1_* \rbrack$, are equal in $H^1(\cF^*)$, as asserted.
\end{proof}

In view of this Theorem, it is natural to ask whether a given cohomology class in $H^1(\cF^*)$ is the infinitesimal of a formal deformation.  This question will be dealt with in Section \ref{sec:integration}.


\section{Extending formal automorphisms and rigidity}
\label{sec:extending}

As in the previous section, $A$ will denote an $S$-algebra with Landweber-Novikov operations $s_*$ and $\cF^* = \cF^*(A)$.

The main purpose of this section is to obtain cohomological conditions under which $A$ is \emph{rigid}, that is, every formal deformation is equivalent to the trivial deformation $s_*$.  To do that, we first have to consider how truncated formal automorphisms can be extended.

\subsection{Formal automorphisms of finite order}
\label{subsec:extending}

Let $m$ be a positive integer.  Inspired by Gerstenhaber-Wilkerson \cite{ger5}, we define a \emph{formal automorphism of order} $m$ on $A$ to be a formal power series
   \[
   \Phi_t ~=~ 1 ~+~ t\phi_1 ~+~ t^2 \phi_2 ~+~ \cdots ~+~ t^m\phi_m
   \]
with $1 = \Id_A$ and each $\phi_i \in \End(A)$, satisfying multiplicativity \eqref{eq:mult} modulo $t^{m+1}$, i.e. the equality \eqref{eq:mult'} holds for $0 \leq n \leq m$ and all $a, b \in A$.  One can think of a formal automorphism as a formal automorphism of order $\infty$.

Given such a $\Phi_t$, we say that it \emph{extends} to a formal automorphism of order $m + 1$ if and only if there exists $\phi_{m + 1} \in \End(A)$ such that the power series
   \begin{equation}
   \label{eq:extension}
   \tilde{\Phi}_t ~=~ \Phi_t ~+~ t^{m+1}\phi_{m+1}
   \end{equation}
is a formal automorphism of order $m + 1$.  Such a $\tilde{\Phi}_t$ is said to be an \emph{extension} of $\Phi_t$ to order $m + 1$.

The question we would like to address here is this: Given a formal automorphism $\Phi_t$ of order $m$, can it be extended to a formal automorphism of order $m + 1$?  It turns out that the obstruction to the existence of such an extension lies in Hochschild cohomology, which we now recall.

Consider $A$ as a unital algebra over the ring of integers $\bZ$, the Hochschild cochain complex $C^*(A,A)$ of $A$ with coefficients in $A$ itself is defined as follows.  The $n$th dimension is the group
   \[
   C^n(A, A) ~=~ \Hom(A^{\otimes n}, A).
   \]
The differential
   \[
   b_{n-1} \colon C^{n-1}(A, A) ~\to~ C^n(A, A)
   \]
is given by the alternating sum
   \begin{multline}
   (b_{n-1}f)(a_1 \otimes \cdots \otimes a_n)
   ~=~ a_1 f(a_2 \otimes \cdots \otimes a_n) ~+~ \\ \sum_{i = 1}^{n-1} (-1)^i f( \cdots a_{i-1} \otimes a_ia_{i+1} \otimes a_{i+2} \cdots) + (-1)^n f(a_1 \otimes \cdots \otimes a_{n-1})a_n.
   \end{multline}
The reader can consult, for example, Weibel \cite{weibel} for more detailed discussions about Hochschild cohomology.  The $n$th Hochschild cohomology group of $A$ over $\bZ$ with coefficients in $A$ itself is denoted by $HH^n(A)$.

We now return to the question of extending a formal automorphism $\Phi_t$ of order $m$ to one of order $m + 1$.  Consider the Hochschild $2$-cochain
   \[
   \Ob(\Phi_t) \colon A^{\otimes 2} ~\to~ A
   \]
given by
   \[
   \Ob(\Phi_t)(a \otimes b) ~=~ - \sum_{i = 1}^m \phi_i(a)\phi_{m+1-i}(b)
   \]
for all $a, b \in A$.  Call $\Ob(\Phi_t)$ the \emph{obstruction class} of $\Phi_t$.

\begin{lemma}
\label{lem:hoch}
The obstruction class $\Ob(\Phi_t)$ is a Hochschild $2$-cocycle.
\end{lemma}

\begin{proof}
We calculate as follows:
   \[
   \begin{split}
   (b_2 \Ob(\Phi_t))&(a \otimes b \otimes c) \\
   &=~ -a\sum_{i=1}^m\phi_i(b)\phi_{m+1-i}(c) ~+~ \sum_{i=1}^m \phi_i(ab)\phi_{m+1-i}(c) \\
   &\qquad \qquad -~ \sum_{i=1}^m \phi_i(a)\phi_{m+1-i}(bc) ~+~ c\sum_{i=1}^m \phi_i(a)\phi_{m+1-i}(b) \\
   &=~ \sum_{\substack{i+j+k ~=~ m+1 \\i,\, k > 0}}\, \phi_i(a)\phi_j(b)\phi_k(c) ~-~ \sum_{\substack{i+j+k ~=~ m+1 \\i,\, k > 0}}\, \phi_i(a)\phi_j(b)\phi_k(c). \\
   &=~ 0
   \end{split}
   \]
Here we used the multiplicativity of $\Phi_t$, namely, \eqref{eq:mult'} for $\phi_i(ab)$ and $\phi_{m+1-i}(bc)$.  Also, $\phi_0 = 1$ as usual.  This shows that the obstruction class of $\Phi_t$ is a Hochschild $2$-cocycle.
\end{proof}

We are now ready to show that the obstruction to extending $\Phi_t$ to a formal automorphism of order $m + 1$ is exactly the cohomology class $\lbrack \Ob(\Phi_t) \rbrack \in HH^2(A)$.

\begin{thm}
\label{thm:ext}
Let $\Phi_t$ be a formal automorphism of order $m$ on $A$.  Then $\Phi_t$ extends to a formal automorphism of order $m + 1$ if and only if the obstruction class $\Ob(\Phi_t)$ is a Hochschild $2$-coboundary.
\end{thm}

\begin{proof}
The existence of an order $m + 1$ extension $\tilde{\Phi}_t$ as in \eqref{eq:extension} is equivalent to the existence of a $\phi_{m+1} \in \End(A)$ for which \eqref{eq:mult'} holds for $n = m + 1$.  This last condition can be rewritten as
   \[
   \Ob(\Phi_t)(a \otimes b)
   ~=~ - \phi_{m+1}(ab) ~+~ a\phi_{m+1}(b) + b\phi_{m+1}(a) ~=~ (b_1 \phi_{m+1})(a \otimes b),
   \]
and $b_1 \phi_{m+1}$ is a Hochschild $2$-coboundary.
\end{proof}

An immediate consequence of this Theorem is that, starting with a derivation $\phi$, the vanishing of the Hochschild cohomology $HH^2(A)$ implies that the formal automorphism $\Phi_t = 1 + t^m \phi$ of order $m \geq 1$ can always be extended to a formal automorphism.  This will be useful when we discuss rigidity below.

\begin{cor}
\label{cor:ext}
Let $m$ be a positive integer and let $\phi$ be a derivation on $A$.  Assume that $HH^2(A) = 0$.  Then there exists a formal automorphism on $A$ of the form
   \[
   \Phi_t ~=~ 1 ~+~ t^m\phi ~+~ t^{m+1}\phi_{m+1} ~+~ t^{m+2}\phi_{m+2} ~+~ \cdots.
   \]
\end{cor}

\begin{proof}
Using the hypothesis that $\phi$ is a derivation, it is straightforward to verify that the formal power series
   \[
   \Phi_t ~=~ 1 ~+~ t^m\phi
   \]
is a formal automorphism of order $m$, i.e.\ it satisfies multiplicativity \eqref{eq:mult'} for $n \leq m$.  As $HH^2(A) = 0$, by Theorem \ref{thm:ext} the obstructions to extending $\Phi_t$ to a formal automorphism vanish, as desired.
\end{proof}

\subsection{Rigidity}
\label{subsec:rigidity}

Following the terminology in \cite{ger1}, an $S$-algebra $A$ is said to be \emph{rigid} if every formal deformation of $A$ is equivalent to the trivial deformation $s_*$.

Using the results above, we will be able to obtain cohomological criterion for the rigidity of $A$.  First we need the following consequence of Corollary \ref{cor:ext}.

\begin{cor}
\label{cor1:rigidity}
Let $k$ be a positive integer and let
   \[
   \sigma^t_* ~=~ s_* ~+~ t^ks^k_* ~+~ t^{k+1}s^{k+1}_* ~+~ \cdots
   \]
be a formal deformation of $A$ in which $s^k_* \in \cF^1$ is a $1$-coboundary.  Suppose that $HH^2(A) = 0$.  Then there exists a formal automorphism of the form
   \begin{equation}
   \label{eq:rigid}
   \Phi_t ~=~ 1 ~-~ t^k \phi_k ~+~ t^{k+1}\phi_{k+1} ~+~ \cdots
   \end{equation}
such that
   \begin{equation}
   \label{eq1:rigid}
   \Phi^{-1}_t \sigma^t_* \Phi_t
   ~\equiv~ s_* \pmod{t^{k+1}}.
   \end{equation}
\end{cor}

\begin{proof}
By the hypothesis on $s^k_*$, there exists a derivation $\phi_k \in \Der(A)$ such that
   \[
   d^0\phi_k ~=~ s_*\phi_k ~-~ \phi_k s_* ~=~ s^k_*.
   \]
Since $-\phi_k$ is also a derivation on $A$, Corollary \ref{cor:ext} implies that there exists a formal automorphism $\Phi_t$ as in \eqref{eq:rigid}.  Computing modulo $t^{k+1}$, we have
   \[
   \begin{split}
   \Phi^{-1}_t \sigma^t_* \Phi_t
   &~\equiv~ (1 ~+~ t^k\phi_k)(s_* ~+~ t^ks^k_*)(1 ~-~ t^k\phi_k) \\
   &~\equiv~ s_* ~+~ t^k(s^k_* ~+~ \phi_k s_* ~-~ s_*\phi_k) \\
   &~\equiv~ s_*.
   \end{split}
   \]
This proves the Corollary.
\end{proof}

By Proposition \ref{prop:aut}, the formal power series $\Phi^{-1}_t \sigma^t_* \Phi_t$ is actually a formal deformation that is equivalent to $\sigma^t_*$.  Therefore, by applying this Corollary repeatedly, we obtain the following cohomological criterion for $A$ to be rigid.

\begin{cor}
\label{cor2:rigid}
If $HH^2(A) = 0$ and $H^1(\cF^*(A)) = 0$, then $A$ is rigid.
\end{cor}

We should point out that in the deformation theory of some other kinds of algebras, rigidity is usually guaranteed by the vanishing of just one cohomology group, usually an $H^1$ or an $H^2$.


\section{Extending cocycles}
\label{sec:integration}

As before, $A$ will denote an arbitrary but fixed $S$-algebra, and $\cF^* = \cF^*(A)$.  The Landweber-Novikov operations on $A$ are given by a function $s_*$ as in \eqref{eq:s}.

Recall from Theorem \ref{thm:inf} that the infinitesimal of a formal deformation is a $1$-cocycle in $\cF^*$.  The purpose of this section is to answer the questions: (1) What is the obstructions to extending a $1$-cocycle to a formal deformation?  (2) What is the obstructions to two such extensions being equivalent?

We will break these questions into a sequence of smaller questions, each of which is dealt with in an obstruction theoretic way.

\subsection{Formal deformations of finite order}
\label{subsec:finite}

First we need some definitions.  Let $m$ be a positive integer.  As in the previous section, a \emph{formal deformation of order} $m$ of $A$ is a formal power series
   \begin{equation}
   \label{eq:def m}
   \sigma^t_* ~=~ s_* ~+~ ts^1_* ~+~ \cdots + t^ms^m_*,
   \end{equation}
in which each $s^i_*$ $(i \geq 1)$ is in $\cF^1$ such that the Cartan formula \eqref{eq:Cartan} and the Product formula \eqref{eq:Product} are satisfied modulo $t^{m+1}$.  In other words, \eqref{eq:Cartan'} and \eqref{eq:Product'} are satisfied for all $n \leq m$.

We say that a formal deformation $\sigma^t_*$ of order $m$ \emph{extends} to a formal deformation of order $M > m$ if and only if there exist $1$-cochains $s^{m+1}_*, \ldots, s^M_*$ such that the power series
   \begin{equation}
   \label{eq:m ext}
   \tilde{\sigma}^t_* ~=~ \sigma^t_* ~+~ t^{m+1}s^{m+1}_* ~+~ \cdots ~+~ t^Ms^M_*
   \end{equation}
is a formal deformation of order $M$.  Call $\tilde{\sigma}^t_*$ an order $M$ extension of $\sigma^t_*$.

One can think of a formal deformation as a formal deformation of order $\infty$.  Given a $1$-cocycle $s^1_*$, the problem of extending the formal deformation
   \[
   \sigma^t_* ~=~ s_* ~+~ ts^1_*
   \]
of order $1$ to a formal deformation can be thought of as extending $\sigma^t_*$ to order $2$, then order $3$, and so forth.  We will do this by identifying the obstruction to extending an order $m$ formal deformation to one of order $m + 1$.

Let $\sigma^t_*$ be a formal deformation of order $m$ as in $\eqref{eq:def m}$.  Consider the $2$-cochain
   \[
   \Ob(\sigma^t_*) ~=~ (\Ob_0(\sigma^t_*), \Ob_1(\sigma^t_*)) \in \cF^2 = \cF^2_0 \times \cF^2_1,
   \]
whose components are defined by the conditions,
   \begin{equation}
   \label{eq:ob0}
   \Ob_0(\sigma^t_*)(\alpha, \beta) ~=~ - \sum_{i = 1}^m s^i_\alpha s^{m+1-i}_\beta \quad (\alpha, \, \beta \in \cE)
   \end{equation}
and
   \begin{equation}
   \label{eq:ob1}
   \Ob_1(\sigma^t_*)_\alpha(a \otimes b) ~=~ - \sum_{i=1}^m \sum_{\beta + \gamma = \alpha} s^i_\beta(a) s^{m+1-i}_\gamma(b) \quad (\alpha \in \cE,\, a, b \in A).
   \end{equation}
We call $\Ob(\sigma^t_*)$ the \emph{obstruction class} of $\sigma^t_*$.

\begin{lemma}
\label{lem1:int}
The obstruction class $\Ob(\sigma^t_*)$ is a $2$-cocycle.
\end{lemma}

\begin{proof}
We need to show that $d^2_i \Ob_i(\sigma^t_*) = 0$ for $i = 0, 1$.  Since $\sigma^t_*$ is the only formal deformation we are dealing with, we abbreviate $\Ob_i(\sigma^t_*)$ to $\Ob_i$.  For $i = 0$, let $\alpha, \beta, \gamma$ be exponential sequences.  Then we have
   \begin{equation}
   \label{eq:2cocycle}
   \begin{split}
   (d^2_0 \Ob_0)&(\alpha, \beta, \gamma)
   ~=~ -s_\alpha\biggl(\sum_{i=1}^m s^i_\beta s^{m+1-i}_\gamma\biggr) ~+~ \sum_{\delta \in P(\alpha, \beta)} n_\delta\biggl(\sum_{i=1}^m s^i_\delta s^{m+1-i}_\gamma\biggr) \\
  &\qquad \qquad \qquad ~-~ \sum_{\varepsilon \in P(\beta, \gamma)} n_{\varepsilon}\biggl(\sum_{i=1}^m s^i_\alpha s^{m+1-i}_\varepsilon\biggr) ~+~ \biggl(\sum_{i=1}^m s^i_\alpha s^{m+1-i}_\beta\biggr)s_\gamma.
   \end{split}
   \end{equation}
This is a sum of four terms.  Using the Product formula \eqref{eq:Product'}, the second term can be rewritten as
   \[
   \sum_{i=1}^m \biggl(\sum_{\delta \in P(\alpha, \beta)} n_\delta s^i_\delta\biggr)s^{m+1-i}_\gamma
   ~=~ \sum_{i=1}^m \sum_{l=0}^i s^l_\alpha s^{i-l}_\beta s^{m+1-i}_\gamma.
   \]
In particular, the $m$ summands in it corresponding to $l = 0$ cancel out with the first term in \eqref{eq:2cocycle}.  Consequently, the sum of the first two terms in \eqref{eq:2cocycle} becomes
   \begin{equation}
   \label{eq2:2cocycle}
   \sum_{\substack{i + j + k ~=~ m + 1 \\ i, k > 0}} s^i_\alpha s^j_\beta s^k_\gamma.
   \end{equation}
An analogous argument applied to the third and fourth terms in \eqref{eq:2cocycle} shows that their sum is exactly \eqref{eq2:2cocycle} with the opposite sign.  It follows that $d^2_0 \Ob_0 = 0$, as desired.

The proof of $d^2_1 \Ob_1 = 0$ is quite similar to the argument above and the proof of Lemma \ref{lem:hoch}.  Indeed, $(d^2_1 \Ob_1)_\alpha(a \otimes b \otimes c)$ (see \eqref{eq:d21}) is again the sum of four terms.  Using the Cartan formula \eqref{eq:Cartan'}, its second term can be rewritten as
   \begin{equation}
   \label{eq3:2cocycle}
   \begin{split}
   -(\Ob_1)_\alpha(ab \otimes c)
   &~=~ \sum_{i=1}^m \sum_{\beta + \gamma = \alpha} s^i_\beta(ab) s^{m+1-i}_\gamma(c) \\
   &~=~ \sum_{i=1}^m \sum_{\beta + \gamma = \alpha}\biggl(\sum_{j=0}^i \sum_{\kappa + \lambda = \beta} s^j_\kappa(a) s^{i-j}_\lambda(b)\biggr) s^{m+1-i}_\gamma(c).
   \end{split}
   \end{equation}
The summands corresponding to $j = 0$ cancel out with the first of the four terms in $(d^2_1 \Ob_1)_\alpha(a \otimes b \otimes c)$.  In particular, the sum of the first two terms in  $(d^2_1 \Ob_1)_\alpha(a \otimes b \otimes c)$ is
   \begin{equation}
   \label{eq4:2cocycle}
   \sum s^i_\beta(a) s^j_\gamma(b) s^k_\delta(c).
   \end{equation}
The summation is taken over all triples $(i, j, k)$ of non-negative integers and all triples $(\beta, \gamma, \delta)$ of exponential sequences for which
   \[
   i ~+~ j ~+~ k ~=~ m ~+~ 1 \qquad (i,\, k > 0)
   \]
and
   \[
   \beta ~+~ \gamma ~+~ \delta ~=~ \alpha.
   \]
A similar argument applies to the last two terms in $(d^2_1 \Ob_1)_\alpha(a \otimes b \otimes c)$, showing that their sum is exactly \eqref{eq4:2cocycle} with the opposite sign.  It follows that $d^2_1 \Ob_1 = 0$, as expected.

This finishes the proof of the Lemma.
\end{proof}

As the obstruction class of $\sigma^t_*$ is a $2$-cocycle, it represents a cohomology class in $H^2(\cF^*)$.

We are now ready to present the main result of this section, which identifies the obstruction to extending a formal deformation of order $m$ to one of order $m + 1$.

\begin{thm}
\label{thm1:int}
Let $\sigma^t_*$ be a formal deformation of order $m$ of $A$.  Then $\sigma^t_*$ extends to a formal deformation of order $m + 1$ if and only if the cohomology class $\lbrack \Ob(\sigma^t_*) \rbrack \in H^2(\cF^*(A))$ vanishes.
\end{thm}

\begin{proof}
Indeed, the existence of an order $m + 1$ extension $\tilde{\sigma}^t_*$ of $\sigma^t_*$, as in \eqref{eq:extension} with $M = m + 1$, is equivalent to the existence of a $1$-cochain $s^{m+1}_* \in \cF^1$ for which the Cartan formula \eqref{eq:Cartan'} and the Product formula \eqref{eq:Product'} both hold for $n = m + 1$.  Simply by rearranging terms, the Cartan formula \eqref{eq:Cartan'} when $n = m + 1$ can be rewritten as
   \[
   (\Ob_1(\sigma^t_*))_\alpha(a \otimes b) ~=~ (d^1_1 s^{m+1}_*)_\alpha(a \otimes b).
   \]
Similarly, the Product formula \eqref{eq:Product'} when $n = m + 1$ is equivalent to
   \[
   (\Ob_0(\sigma^t_*))(\alpha, \beta) ~=~ (d^1_0 s^{m+1}_*)(\alpha, \beta).
   \]
These two conditions together are equivalent to
   \begin{equation}
   \label{eq:cob}
   \Ob(\sigma^t_*) ~=~ d^1 s^{m+1}_*,
   \end{equation}
i.e., the obstruction class is a $2$-coboundary.  Since $\Ob(\sigma^t_*)$ is a $2$-cocycle (Lemma \ref{lem1:int}), the Theorem is proved.
\end{proof}

Applying this Theorem repeatedly, we obtain the obstructions to extending a $1$-cocycle to a formal deformation.

\begin{cor}
\label{cor1:int}
Let $s^1_* \in \cF^1(A)$ be a $1$-cocycle.  Then there exist a sequence of classes $\omega_1, \omega_2, \ldots \in H^2(\cF^*(A))$ for which $\omega_n$ $(n > 1)$ is defined if and only if $\omega_1, \ldots , \omega_{n-1}$ are defined and equal to $0$.  Moreover, the formal deformation
   \[
   \sigma^t_* ~=~ s_* ~+~ ts^1_*
   \]
of order $1$ on $A$ extends to a formal deformation if and only if $\omega_i$ is defined and equal to $0$ for each $i = 1, 2, \ldots$.
\end{cor}

Since these obstructions lie in $H^2(\cF^*(A))$, the triviality of this group implies that extensions always exist.

\begin{cor}
\label{cor2:int}
If $H^2(\cF^*(A)) = 0$, then every formal deformation of order $m \geq 1$ of $A$ extends to a formal deformation.
\end{cor}

\subsection{Equivalence of formal deformations of finite order}
\label{subsec:equiv}

Let $\sigma^t_*$ be a formal deformation of $A$ of order $m \geq 1$, and let $\tilde{\sigma}^t_*$ and $\bar{\sigma}^t_*$ be two order $m + 1$ extensions of it.  Are the two extensions equivalent?  This is the question that we would like to address in this final section.

First, we need a definition of equivalence.  Two formal deformation $\sigma^t_*$ and $\tilde{\sigma}^t_*$ of order $m$ of $A$ are said to be \emph{equivalent} if and only if there exists a formal automorphism $\Phi_t$ of order $m$ such that the equality \eqref{eq:equivalence} holds modulo $t^{m+1}$.

This is a well-defined equivalence relation on the set of formal deformations of order $m$ of $A$.  In fact, it is easy to see that a formal automorphism $\Phi_t$ of order $m$ has a formal inverse $\Phi^{-1}_t$ as in \eqref{eq:inv} which is also a formal automorphism of order $m$.

Now let $\sigma^t_* = s_* + \cdots + t^ms^m_*$ be a formal deformation of order $m$ of $A$, and let $\tilde{\sigma}^t_* = \sigma^t_* + t^{m+1}\tilde{s}^{m+1}_*$ and $\bar{\sigma}^t_* = \sigma^t_* + t^{m+1}\bar{s}^{m+1}_*$ be two order $m + 1$ extensions of $\sigma^t_*$.  Then \eqref{eq:cob} in the proof of Theorem \ref{thm1:int} tells us that
   \[
   d^1 \tilde{s}^{m+1}_* ~=~ \Ob(\sigma^t_*) ~=~ d^1 \bar{s}^{m+1}_*.
   \]
It follows that $(\tilde{s}^{m+1}_* - \bar{s}^{m+1}_*) \in \cF^1$ is a $1$-cocycle, and it makes sense to consider the cohomology class in $H^1(\cF^*)$ represented by it.

\begin{thm}
\label{thm2:int}
If the cohomology class $\lbrack \tilde{s}^{m+1}_* - \bar{s}^{m+1}_* \rbrack \in H^1(\cF^*)$ vanishes, then the formal deformations, $\tilde{\sigma}^t_*$ and $\bar{\sigma}^t_*$, of order $m + 1$ are equivalent.
\end{thm}

\begin{proof}
The argument is quite similar to the proofs of Theorem \ref{thm:inf} and Corollary \ref{cor1:rigidity}.  In fact, by hypothesis, there exists a derivation $\phi$ on $A$ such that
   \[
   \tilde{s}^{m+1}_* ~-~ \bar{s}^{m+1}_* ~=~ d^0\phi ~=~ s_*\phi ~-~ \phi s_*.
   \]
We have the formal automorphism
   \[
   \Phi_t ~=~ 1 ~+~ t^{m+1}\phi
   \]
of order $m + 1$ on $A$.  Computing modulo $t^{m+2}$, we have
   \[
   \begin{split}
   \Phi^{-1}_t \bar{\sigma}^t_* \Phi_t
   &~\equiv~ (1 ~-~ t^{m+1}\phi)(\sigma^t_* ~+~ t^{m+1}\bar{s}^{m+1}_*)(1 ~+~ t^{m+1}\phi) \\
   &~\equiv~ \sigma^t_* ~+~ t^{m+1}(\bar{s}^{m+1}_* ~+~ s_*\phi ~-~ \phi s_*) \\
   &~\equiv~ \sigma^t_* ~+~ t^{m+1}\tilde{s}^{m+1}_* \\
   &~\equiv~ \tilde{\sigma}^t_*,
   \end{split}
   \]
as desired.
\end{proof}

This obstruction theoretic result is less satisfactory than the analogous Theorem \ref{thm:ext} and Theorem \ref{thm1:int} in that the author is not sure whether the equivalence of $\tilde{\sigma}^t_*$ and $\bar{\sigma}^t_*$ would imply the vanishing of the class $\lbrack \tilde{s}^{m+1}_* - \bar{s}^{m+1}_* \rbrack$.

\section*{Acknowledgement}
The author thanks the referee for reading an earlier version of this paper.


\end{document}